# A Family of Estimators of Population Mean Using Information on Point Bi-Serial and Phi Correlation Coefficient

**Sachin Malik** and **Rajesh Singh***

Department of Statistics, Banaras Hindu University
Varanasi-221005, India
*Corresponding Author
(sachinkurava999@gmail.com, rsinghstat@gmail.com)

**ABSTRACT**

*This paper deals with the problem of estimating the finite population mean when some information on two auxiliary attributes are available. It is shown that the proposed estimator is more efficient than the usual mean estimator and other existing estimators. The study is also extended to two-phase sampling. The results have been illustrated numerically by taking empirical population considered in the literature.*

## 1. INTRODUCTION

There are some situations when in place of one auxiliary attribute, we have information on two qualitative variables. For illustration, to estimate the hourly wages we can use the information on marital status and region of residence (see Gujrati and Sangeetha (2007), page-311). Here we assume that both auxiliary attributes have significant point bi-serial correlation with the study variable and there is significant phi-correlation (see Yule (1912)) between the auxiliary attributes. The use of auxiliary information can increase the precision of an estimator when study variable Y is highly correlated with auxiliary variables X. Naik and Gupta (1996) introduced a ratio estimator when the study variable and the auxiliary attribute are positively correlated. Jhajj et al. (2006) suggested a family of estimators for the population mean in single and two phase sampling when the study variable and auxiliary attribute are positively correlated. Shabbir and Gupta (2007), Singh et al. (2008), Singh et al. (2010) and Abd-Elfattah et al. (2010) have considered the problem of estimating population mean $\overline{Y}$ taking into consideration the point biserial correlation between auxiliary attribute and study variable.

In order to have an estimate of the study variable y, assuming the knowledge of the population proportion P, Naik and Gupta (1996) and Singh et al. (2007) respectively proposed following estimators

$$t_1 = \bar{y}\left(\frac{P_1}{p_1}\right) \tag{1.1}$$

$$t_2 = \bar{y}\left(\frac{p_2}{P_2}\right) \tag{1.2}$$

$$t_3 = \bar{y}\exp\left(\frac{P_1 - p_1}{P_1 + p_1}\right) \tag{1.3}$$

$$t_4 = \bar{y}\exp\left(\frac{p_2 - P_2}{p_2 + P_2}\right) \tag{1.4}$$

The bias and MSE expression's of the estimator's $t_i$ (i=1, 2, 3, 4) up to the first order of approximation are, respectively, given by

$$B(t_1) = \bar{Y} f_1 C_{p_1}^2 \left[1 - K_{pb_1}\right] \tag{1.5}$$

$$B(t_2) = \bar{Y} f_1 K_{pb_2} C_{p_2}^2 \tag{1.6}$$

$$B(t_3) = \bar{Y} f_1 \frac{C_{p_2}^2}{2}\left[\frac{1}{4} - K_{pb_2}\right] \tag{1.7}$$

$$B(t_4) = \bar{Y} f_1 \frac{C_{p_2}^2}{2}\left[\frac{1}{4} + K_{pb_2}\right] \tag{1.8}$$

$$\text{MSE}(t_1) = \bar{Y}^2 f_1\left[C_y^2 + C_{p_1}^2\left(1 - 2K_{pb_1}\right)\right] \tag{1.9}$$

$$\text{MSE}(t_2) = \bar{Y}^2 f_1\left[C_y^2 + C_{p_1}^2\left(1 + 2K_{pb_2}\right)\right] \tag{1.10}$$

$$\text{MSE}(t_3) = \bar{Y}^2 f_1\left[C_y^2 + C_{p_1}^2\left(\frac{1}{4} - K_{pb_2}\right)\right] \tag{1.11}$$

$$\text{MSE}\left(t_4\right)=\overline{Y}^2 f_1\left[C_y^2+C_{p_2}^2\left(\frac{1}{4}+K_{pb_2}\right)\right] \tag{1.12}$$

where, $f_1=\frac{1}{n}-\frac{1}{N}$ , $S_{\phi_j}^2=\frac{1}{N-1}\sum_{i=1}^{N}\left(\phi_{ji}-P_j\right)^2$, $S_{y\phi_j}=\frac{1}{N-1}\sum_{i=1}^{N}\left(y_i-\overline{Y}\right)\left(\phi_{ji}-P_j\right)$,

$$\rho_{pb_j}=\frac{S_{y\phi_j}}{S_y S_{\phi_j}},\ C_y=\frac{S_y}{\overline{Y}},\ C_{pj}=\frac{S_{\phi_j}}{P_j};\ (j=1,2),$$

$$K_{pb_1}=\rho_{pb_1}\frac{C_y}{C_{p_1}},\ K_{pb_2}=\rho_{pb_2}\frac{C_y}{C_{p_2}}.$$

$s_{\phi_1\phi_2}=\frac{1}{n-1}\sum_{i=1}^{n}\left(\phi_{1i}-p_1\right)\left(\phi_{2i}-p_2\right)$ and $\rho_\phi=\frac{s_{\phi_1\phi_2}}{s_{\phi_1}s_{\phi_2}}$ be the sample phi-covariance and phi-correlation between $\phi_1$ and $\phi_2$ respectively, corresponding to the population phi-covariance and phi-correlation $S_{\phi_1\phi_2}=\frac{1}{N-1}\sum_{i=1}^{N}\left(\phi_{1i}-P_1\right)\left(\phi_{2i}-P_2\right)$

$$\text{and } \rho_\phi=\frac{S_{\phi_1\phi_2}}{S_{\phi_1}S_{\phi_2}}$$

Following Naik and Gupta (1996) and Singh et al. (2007), we propose the estimators $t_5$ and $t_6$ as

$$t_5=\overline{y}\left(\frac{P_1}{p_1}\right)^{\alpha_1}\left(\frac{P_2}{p_2}\right)^{\alpha_2} \tag{1.13}$$

$$t_6=\overline{y}\exp\left(\frac{P_1-p_1}{P_1+p_1}\right)^{\beta_1}\exp\left(\frac{p_2-P_2}{p_2+P_2}\right)^{\beta_2} \tag{1.14}$$

where $\alpha_1,\alpha_2,\beta_1$ and $\beta_2$ are real constants.

## 2. BIAS AND MSE of $t_5$ and $t_6$

To obtain the bias and MSE of the estimators $t_5$ and $t_6$ to the first degree of approximation, we define

$$e_0=\frac{\overline{y}-\overline{Y}}{\overline{Y}},\ e_1=\frac{p_1-P_1}{P_1},\ e_2=\frac{p_2-P_2}{P_2}$$

such that, $E(e_i) = 0$; $i = 0, 1, 2.$

Also,

$$E(e_0^2) = f_1 C_y^2,\ E(e_1^2) = f_1 C_{p_1}^2,\ E(e_2^2) = f_1 C_{p_2}^2,$$

$$E(e_0 e_1) = f_1 K_{pb_1} C_{p_1}^2,\ E(e_0 e_2) = f_1 K_{pb_2} C_{p_2}^2,\quad E(e_1 e_2) = f_1 K_{\phi} C_{p_2}^2,$$

$$K_{pb_1} = \rho_{pb_1} \frac{C_y}{C_{p_1}},\ K_{pb_2} = \rho_{pb_2} \frac{C_y}{C_{p_2}},\ K_{\phi} = \rho_{\phi} \frac{C_{p_1}}{C_{p_2}}$$

Expressing equation (1.13) in terms of e's, we have

$$t_5 = \overline{Y}\left[(1+e_0)(1+e_1)^{-\alpha_1}(1+e_2)^{-\alpha_2}\right]$$

$$= \overline{Y}\left[1 - \alpha_1 e_1 - \alpha_2 e_2 + \frac{\alpha_1(\alpha_1+1)}{2} e_1^2 + \frac{\alpha_2(\alpha_2+1)}{2} e_2^2 + \alpha_1\alpha_2 e_1 e_2 - \alpha_1 e_0 e_1 - \alpha_2 e_0 e_2\right] \quad (2.1)$$

Subtracting $\overline{Y}$ from both the sides of equation (2.1) and then taking expectation of both sides, we get the bias of the estimator $t_5$, up to the first order of approximation, as

$$B(t_5) = \overline{Y} f_1 \left( C_{p_1}^2 \left[ \frac{\alpha_1^2}{2} + \frac{\alpha_1}{2} - \alpha_1 k_{pb_1} \right] + C_{p_2}^2 \left[ \frac{\alpha_2^2}{2} + \frac{\alpha_2}{2} - \alpha_2 k_{pb_2} + \alpha_1 \alpha_2 k_{\phi} \right] \right) \quad (2.2)$$

From (2.1), we have

$$(t_5 - \overline{Y}) = \overline{Y}\left[e_0 - \alpha_1 e_1 - \alpha_2 e_2\right] \quad (2.3)$$

Squaring both sides of (2.3) and then taking expectation, we get, MSE of the estimator $t_5$, up to the first order of approximation, as

$$MSE(t_5) = \overline{Y}^2 f_1 \left[ C_y^2 + C_{p_1}^2 \left( \alpha_1^2 - 2\alpha_1 K_{pb_1} \right) + C_{p_2}^2 \left( \alpha_2^2 - 2\alpha_2 K_{pb_2} + 2\alpha_1\alpha_2 K_{\phi} \right) \right] \quad (2.4)$$

To obtain the bias and MSE of $t_6$, to the first order of approximation, we express equation (1.14) in term of e's, as

$$t_6 = \overline{Y}\left(1 + e_0 \exp\left(\frac{-e_1}{2+e_1}\right)^{\beta_1} \exp\left(\frac{e_2}{2+e_2}\right)^{\beta_2}\right)$$

$$= \overline{Y}\left[1 + e_0 - \frac{\beta_1 e_1}{2} + \frac{\beta_1^2 e_1^2}{4} - \frac{\beta_2 e_2}{2} - \frac{\beta_1\beta_2 e_1 e_2}{4} + \frac{\beta_2^2 e_2^2}{4} + \frac{\beta_2 e_0 e_2}{2} - \frac{\beta_1 e_0 e_1}{2}\right] \tag{2.5}$$

Subtracting $\overline{Y}$ from both sides of equation (2.5) and then taking expectation of both sides, we get the bias of the estimator $t_6$ up to the first order of approximation, as

$$B(t_6) = \overline{Y} f_1\left[C_{p_1}^2\left(\frac{\beta_1^2}{4} - \frac{\beta_1}{2}K_{pb_1}\right) + C_{p_2}^2\left(\frac{\beta_2^2}{4} + \frac{\beta_2}{2}K_{pb_2} - \frac{\beta_1\beta_2}{4}K_\phi\right)\right] \tag{2.6}$$

From (2.5), we have

$$(t_6 - \overline{Y}) = \overline{Y}\left[e_0 - \frac{\beta_1 e_1}{2} - \frac{\beta_2 e_2}{2}\right] \tag{2.7}$$

Squaring both sides of (2.7) and then taking expectation, we get the MSE of the estimator $t_6$ up to the first order of approximation, as

$$MSE(t_6) = \overline{Y}^2 f_1\left[C_y^2 + C_{p_1}^2\left(\frac{\beta_1^2}{4} - \beta_1 K_{pb_1}\right) + C_{p_2}^2\left(\frac{\beta_2^2}{4} - \frac{\beta_1\beta_2}{2}K_\phi + \beta_2 K_{pb_1}\right)\right] \tag{2.8}$$

## 3. ANOTHER ESTIMATOR

Following Naik and Gupta (1996) and Singh et al. (2007), we propose another improved estimator $t_p$ as-

$$t_p = w_0\bar{y} + w_1\bar{y}\left(\frac{P_1}{p_1}\right)^{\alpha_1}\left(\frac{P_2}{p_2}\right)^{\alpha_2} + w_2\bar{y}\exp\left(\frac{P_1 - p_1}{P_1 + p_1}\right)^{\beta_1}\exp\left(\frac{p_2 - P_2}{p_2 + P_2}\right)^{\beta_2} \quad (3.1)$$

where $\alpha_1, \alpha_2$, $\beta_1$ and $\beta_2$ are real constants and $w_i (i = 0,1,2)$ are suitably chosen constants whose values are to be determined later.

Expressing (3.1) in terms of e's, we have

$$t_p = \bar{Y}(1+e_0)\left[w_0 + w_1(1+e_1)^{-\alpha_1}(1+e_2)^{-\alpha_2} + w_2\exp\left(\frac{-\beta_1 e_1}{2}\right)\exp\left(\frac{\beta_2 e_2}{2}\right)\right] \quad (3.2)$$

Expanding the right hand side of equation (3.2) and retaining terms, up to second power of e's, we have

$$t_p = \bar{Y}\left[1 + e_0 + w_1\left(\frac{\alpha_1(1+\alpha_1)}{2}e_1^2 + \frac{\alpha_2(1+\alpha_2)}{2}e_2^2 - \alpha_1 e_1 - \alpha_2 e_2 + \alpha_1\alpha_2 e_1 e_2 - \alpha_1 e_0 e_1 - \alpha_2 e_0 e_2\right)\right.$$

$$\left. + w_2\left(\frac{\beta_1^2 e_1^2}{4} - \frac{\beta_1 e_1}{2} - \frac{\beta_1 e_0 e_1}{2} - \frac{\beta_1\beta_2 e_1 e_2}{4} + \frac{\beta_2^2 e_2^2}{4} + \frac{\beta_2 e_2}{2} + \frac{\beta_2 e_0 e_2}{2}\right)\right] \quad (3.3)$$

Subtracting $\bar{Y}$ from both sides of (3.3) and then taking expectation of both sides, we get the bias of the estimator $t_p$, up to the first order of approximation as

$$B(t_p) = \bar{Y}f_1\left\{w_1\left(C_{p_1}^2\left[\frac{\alpha_1^2}{2} + \frac{\alpha_1}{2} - \alpha_1 k_{pb_1}\right] + C_{p_2}^2\left[\frac{\alpha_2^2}{2} + \frac{\alpha_2}{2} - \alpha_2 k_{pb_2} + \alpha_1\alpha_2 k_\phi\right]\right)\right.$$

$$\left. + w_2\left(C_{p_1}^2\left[\frac{\beta_1^2}{4} - \frac{\beta_1}{2}k_{pb_1}\right] + C_{p_2}^2\left[\frac{\beta_2^2}{4} + \frac{\beta_2}{2}k_{pb_2} - \frac{\beta_1\beta_2}{4}k_\phi\right]\right)\right\} \quad (3.4)$$

From (3.3), we have

$$(t_p - \bar{Y}) = \bar{Y}\left[e_0 + w_1(-\alpha_1 e_1 - \alpha_2 e_2) + w_2\left(-\frac{\beta_1 e_1}{2} + \frac{\beta_2 e_2}{2}\right)\right] \quad (3.5)$$

Squaring both sides of (3.6) and then taking expectation, we get MSE of the estimator $t_p$ up to the first order of approximation, as

$$MSE\left(t_p\right)=\overline{Y}^2 f_1\left[C_y^2+w_1^2A_1+w_2^2A_2-2w_1A_3-w_2A_4+w_1w_2A_5\right] \tag{3.6}$$

where,

$$\left.\begin{aligned} w_1 &= \frac{4A_2A_3-A_4A_5}{4A_1A_2-A_5^2} \\ w_2 &= \frac{4A_2A_3-A_4A_5}{4A_1A_2-A_5^2} \end{aligned}\right\} \tag{3.7}$$

and

$$\left.\begin{aligned} A_1 &= \alpha_1^2C_{p_1}^2+\alpha_2^2C_{p_2}^2+2\alpha_1\alpha_2k_\phi C_{p_2}^2 \\ A_2 &= \frac{1}{4}\left[\beta_1^2c_{p_1}^2+\beta_2^2C_{p_2}^2-2\beta_1\beta_2k_\phi C_{p_2}^2\right] \\ A_3 &= \alpha_1k_{pb_1}C_{p_1}^2+\alpha_2k_{pb_2}C_{p_2}^2 \\ A_4 &= \beta_1k_{pb_1}C_{p_1}^2-\beta_2k_{pb_2}C_{p_2}^2 \\ A_5 &= \alpha_1\beta_1C_{p_1}^2-\alpha_2\beta_2C_{p_2}^2+\alpha_2\beta_1k_\phi C_{p_2}^2-\alpha_1\beta_2k_\phi C_{p_2}^2 \end{aligned}\right\} \tag{3.8}$$

## 4. EMPIRICAL STUDY

Data: (Source: Government of Pakistan (2004))

The population consists rice cultivation areas in 73 districts of Pakistan. The variables are defined as:

Y= rice production (in 000' tonnes, with one tonne = 0.984 ton) during 2003,

$P_1$ = production of farms where rice production is more than 20 tonnes during the year 2002, and

$P_2$ = proportion of farms with rice cultivation area more than 20 hectares during the year 2003.

For this data, we have

N=73, $\overline{Y}$ =61.3, $P_1$=0.4247, $P_2$=0.3425, $S_y^2$=12371.4, $S_{\phi_1}^2$ =0.225490, $S_{\phi_2}^2$ =0.228311,

$\rho_{pb_1}$ =0.621, $\rho_{pb_2}$ =0.673, $\rho_\phi$=0.889.

**Table 4.1: PRE of different estimators of $\overline{Y}$ with respect to $\overline{y}$**

| Choice of scalars | | | | | | | | | |
|---|---|---|---|---|---|---|---|---|---|
| $w_0$ | $w_1$ | $w_2$ | $\alpha_1$ | $\alpha_2$ | $\beta_1$ | $\beta_2$ | **Estimator** | **MSE** | **PRE'S** |
| 1 | 0 | 0 | | | | | $\overline{y}$ | 655.28 | 100.00 |
| 0 | 1 | 0 | 1 | 0 | | | $t_1$ | 402.80 | 162.68 |
| | | | 0 | 1 | | | $t_2$ | 1392.16 | 47.66 |
| | | | -1 | 1 | | | $t_5$ | 580.01 | 112.97 |
| 0 | 0 | 1 | | | 1 | 0 | $t_3$ | 462.07 | 141.80 |
| | | | | | 0 | 1 | $t_4$ | 1091.20 | 60.05 |
| | | | | | 1 | -1 | $t_6$ | 363.03 | 180.50 |
| $w_0$ | $w_1$ | $w_2$ | 1 | 1 | 1 | 1 | $t_p$ | 356.87 | 183.60 |

## 5. DOUBLE SAMPLING

It is assumed that the population proportion $P_1$ for the first auxiliary attribute $\phi_1$ is unknown but the same is known for the second auxiliary attribute $\phi_2$. When $P_1$ is unknown, it is some times estimated from a preliminary large sample of size $n'$ on which only the attribute $\phi_1$ is measured. Then a second phase sample of size n (n< $n'$) is drawn and Y is observed.

$$\text{Let } p'_j = \frac{1}{n}\sum_{i=1}^{n'}\phi_{ji}, (j = 1,2).$$

The estimator's $t_1$, $t_2$, $t_3$ and $t_4$ in two-phase sampling take the following form

$$t_{d1} = \overline{y}\left(\frac{p_1'}{p_1}\right) \tag{5.1}$$

$$t_{d2} = \bar{y}\left(\frac{P_2}{p_2^{'}}\right) \quad (5.2)$$

$$t_{d3} = \bar{y}\exp\left(\frac{p_1^{'} - p_1}{p_1^{'} + p_1}\right) \quad (5.3)$$

$$t_{d4} = \bar{y}\exp\left(\frac{p_2^{'} - P_2}{p_2^{'} + P_2}\right) \quad (5.4)$$

The bias and MSE expressions of the estimators $t_{d1}$, $t_{d2}$, $t_{d3}$ and $t_{d4}$ up to first order of approximation, are respectively given as

$$B(t_{d1}) = \bar{Y}f_3C_{p_1}^2\left[1 - k_{pb_1}\right] \quad (5.5)$$

$$B(t_{d2}) = \bar{Y}f_2C_{p_2}^2\left[1 - K_{pb_2}\right] \quad (5.6)$$

$$B(t_{d3}) = \bar{Y}f_3\frac{C_{p_2}^2}{4}\left[1 - K_{pb_2}\right] \quad (5.7)$$

$$B(t_{d4}) = \bar{Y}f_3\frac{C_{p2}^2}{4}\left[1 + K_{pb_2}\right] \quad (5.8)$$

$$\text{MSE}(t_{d1}) = \bar{Y}^2\left[f_1C_y^2 + f_3C_{P_1}^2\left(1 - 2K_{pb_1}\right)\right] \quad (5.9)$$

$$\text{MSE}(t_{d2}) = \bar{Y}^2\left[f_1C_y^2 + f_2C_{p_2}^2\left(1 - 2K_{kp_2}\right)\right] \quad (5.10)$$

$$\text{MSE}(t_{d3}) = \bar{Y}^2\left[f_1C_y^2 + f_3\frac{C_{p_1}^2}{4}\left(1 - 4K_{pb_1}\right)\right] \quad (5.11)$$

$$\text{MSE}(t_{d4}) = \bar{Y}^2\left[f_1C_y^2 + f_3\frac{C_{p_1}^2}{4}\left(1 + 4K_{pb_1}\right)\right] \quad (5.12)$$

where,

$$S^2_{\phi_J} = \frac{1}{n-1}\sum_{i=1}^{n}\left(\phi_{ji} - p_j\right)^2, \quad S'^{\,2}_{\phi_j} = \frac{1}{n'-1}\sum_{i=1}^{n!}\left(\phi_{ji} - p'_j\right)^2,$$

$$f_2 = \frac{1}{n'} - \frac{1}{N}, \quad f_3 = \frac{1}{n} - \frac{1}{n'}.$$

The estimator's $t_5$ and $t_6$, in two phase sampling, takes the following form

$$t_{d5} = \bar{y}\left(\frac{p'_1}{p_1}\right)^{m_1}\left(\frac{P_2}{p'_2}\right)^{m_2} \qquad (5.13)$$

$$t_{d6} = \bar{y}\exp\left(\frac{p'_1 - p_1}{p'_1 + p_1}\right)^{n_1}\exp\left(\frac{p'_2 - P_2}{p'_2 + P_2}\right)^{n_2} \qquad (5.14)$$

Where $m_1, m_2, n_1$ and $n_2$ are real constants.

## 6. BIAS AND MSE OF $t_{d5}$ and $t_{d6}$

To obtain the bias and MSE of $t_{d5}$ and $t_{d6}$ up to first degree of approximation, we define

$$e_o = \frac{\bar{y} - \bar{Y}}{\bar{Y}}, \quad e'_1 = \frac{p'_1 - P_1}{P_1}, \quad e'_2 = \frac{p'_2 - P_2}{P_2}$$

Such that, $E(e_o) = E(e'_1) = E(e'_2) = o.$

Also,

$$E(e_o^2) = f_1 C_y^2, \ E(e_1'^{\,2}) = f_2 C_{p_1}^2, \ E\left(e_2'^{\,2}\right) = f_2 C_{p_2}^2,$$

$$E(e'_1 e'_2) = f_2 C_{p_2}^2, \ E(e_o e'_1) = f_2 K_{pb_1} C_{p_1}^2, \ E(e_o e'_2) = f_2 K_{pb_2} C_{p_2}^2.$$

Expressing equation (5.13) in terms of e's, we have

$$t_{d6} = \bar{Y}(1 + e_o)\left(1 + e'_1\right)^{m_1}(1 + e_1)^{-m_1}\left(1 + e'_2\right)^{-m_2}$$

Expanding the right hand side of above equation and retaining terms up to second power of e's, we have

$$t_{d_5} = \overline{Y}\Big[1 + e_o - m_1 e_1 + \frac{m_1(m_1+1)}{2} e_1^2 - m_1^2 e_1 e_2^{'} + \frac{m_1(m_1-1)}{2} e_1^{'\,2}$$
$$+ m_1 e_1^{'} - m_2 e_2^{'} + \frac{m_2(m_2+1)}{2} e_2^{'\,2} - m_1 e_1 e_0 + m_1 e_0 e_1^{'} - m_2 e_0 e_2^{'}\Big] \quad (6.1)$$

Subtracting $\overline{Y}$ from both sides of (6.1) and then taking expectation, we get the bias of the estimator $t_{d5}$, up to the first order of approximation, as

$$B(t_{d5}) = \overline{Y}\left[f_3 C_{p_1}^2\left(\frac{m_1^2}{2} + \frac{m_1}{2} - m_1 K_{pb_1}\right) + f_2 C_{P_2}^2\left(\frac{m_2^2}{2} + \frac{m_2}{2} - m_2 k_{pb_2}\right)\right] \quad (6.2)$$

From (5.1), we have

$$(t_{d5} - \overline{Y}) = \overline{Y}\left[e_0 - e_1 m_1 + m_1 e_1^{'} - m_2 e_2^{'}\right] \quad (6.3)$$

Squaring both sides of (6.3) and then taking expectations, we get MSE of the estimator $t_{d5}$,

up to the first order of approximation, as

$$MSE(t_{d5}) = \overline{Y}\left[f_1 C_y^2 + f_3 C_{p_1}^2\left(m_1^2 - 2m_1 K_{pb_1}\right) + f_2 C_{p_2 2}^2\left(m_2^2 - 2m_2 K_{pb_2}\right)\right] \quad (6.4)$$

Now to obtain the bias and MSE of $t_{d6}$ to the first order of approximation, we express equation (5.14) in terms of e's

$$t_{d6} = \overline{Y}(1 + e_0)\exp\left(\frac{n_1 e_1^{'}}{2}\right)\exp\left(\frac{-n_1 e_1}{2}\right)\exp\left(\frac{n_2 e_2^{'}}{2}\right)$$

Expanding the right hand side of above equation and retaining terms up to second power of e's, we have

$$t_{d6} = \overline{Y}\left(1 + e_0 + \frac{n_1 e_1^{'}}{2} - \frac{n_1 e_1}{2} + \frac{n_1^2 e_1^2}{4} + \frac{n_2 e_2^{'}}{2} + \frac{n_2 e_2^{\prime 2}}{4} + \frac{n_1 e_0 e_1^{\prime}}{2} - \frac{n_1 e_0 e_1}{2} + \frac{n_2 e_0 e_2^{\prime}}{2}\right) \quad (6.5)$$

Subtracting $\overline{Y}$ from both sides of (6.5) and then taking expectations, we get the bias of the estimator $t_{d6}$ up to the first order of approximation, as

$$B(t_{d6}) = \overline{Y}\left[f_3\left(\frac{n_1^2}{8} + \frac{n_1}{8} - \frac{n_1}{2} K_{pb_1}\right) C_{p_1}^2 + f_2\left(\frac{n_2^2}{8} + \frac{n_2}{8} + \frac{n_2}{2} K_{pb_2}\right)\right] \quad (6.6)$$

From (6.5), we have

$$(t_{d6} - \overline{Y}) = \overline{Y}\left(e_0 + \frac{n_1 e_1'}{2} - \frac{n_1 e_1}{2} + \frac{n_2 e_2'}{2}\right) \tag{6.7}$$

Squaring both sides of (6.7) and then taking expectations, we get the MSE of $t_{d6}$ up to the first order of approximation, as

$$MSE(t_{d6}) = \overline{Y}^2\left[f_1 C_y^2 + f_3\left(\frac{n_1^2}{4} - n_1 K_{pb_1}\right)C_{p_1}^2 + f_2\left(\frac{n_2^2}{4} + n_2 K_{pb_2}\right)C_{p_2}^2\right] \tag{6.8}$$

### 7. IMPROVED ESTIMATOR $t_p$ IN TWO-PHASE SAMPLING

The estimator $t_p$ in double sampling is written as

$$t_{pd} = h_0 \overline{y} + h_1 \overline{y}\left(\frac{p_1'}{p_1}\right)^{m_1}\left(\frac{P_2}{p_2'}\right)^{m_2} + h_2 \exp\left(\frac{p_1' - p_1}{p_1' + p_1}\right)^{n_1} \exp\left(\frac{p_2' - P_2}{p_2' + P_2}\right)^{n_2} \tag{7.1}$$

where, $m_1, m_2, n_1$ and $n_2$ are real constants and $h_i\,(i = 0,1,2)$ are suitably chosen constants whose values are to be determined later.

Expressing (7.1) in terms of e's, we have

$$t_{pd} = h_0 \overline{y} + h_1 \overline{y}(1 + e_1')^{m_1}(1 + e_1)^{-m_2}(1 + e_2')^{-m_2}$$
$$+ h_2 \overline{y} \exp\left(\frac{n_1 e_1'}{2}\right)\exp\left(\frac{-n_1 e_1}{2}\right)\exp\left(\frac{n_2 e_2'}{2}\right) \tag{7.2}$$

Expanding the right hand side of (7.2) and retaining terms up to second power of e's as

$$t_{pd} = \overline{Y}\Bigg[1 + e_0 + h_1\Bigg(\frac{m_1(m_1+1)}{2}e_1^2 - m_1 e_1 - m_1 e_0 e_1 - m_2^2 e_1 e_1' + \frac{m_1(m_1-1)2}{2}e_1'^2 + m_1 e_1'$$
$$+ m_1 e_0 e_1' - m_2 e_2' - m_2 e_0 e_2' + \frac{m_2(m_2+1)}{2}e_2'^2\Bigg)$$
$$+ w_2\left(\frac{n_1 e_1'}{2} - \frac{n_1 e_1}{2} + \frac{n_1^2 e_1^2}{4} + \frac{n_2 e_2'}{2} + \frac{n_2 e_2'^2}{4} + \frac{n_1 e_0 e_1'}{2} - \frac{n_1 e_0 e_1}{2} + \frac{n_2 e_0 e_2'}{2}\right)\Bigg] \tag{7.3}$$

Subtracting $\overline{Y}$ from both the sides of (7.3) then taking expectations on both the sides, we get the bias of the estimator $t_d$ up to the first order of approximation as

$$B(t_{pd}) = \overline{Y}\left[h_1 f_3 C_{p_1}^2\left(\frac{m_1{}^2}{2} + \frac{m_1}{2} + -m_1 K_{pb_1}\right) + h_1 f_2 C_{p_2}^2\left(\frac{m_2^2}{2} + \frac{m_2}{2} - m_2 k_{pb_2}\right)\right.$$
$$\left. + h_3 f_3 C_{p_1}^2\left(\frac{n_1^2}{8} + \frac{n_1}{8} - \frac{n_1}{2}K_{pb_1}\right) + h_3 f_2 C_{p_2}^2\left(\frac{n_2^2}{8} + \frac{n_2}{8} - \frac{n_2}{2}k_{pb_2}\right)\right] \qquad (7.4)$$

From (7.3), we have

$$t_{pd} - \overline{Y} = \overline{Y}\left[e_0 + h_1\left(-m_1 e_1 + m_1 e_1' - m_1 e_2'\right) + h_2\left(\frac{n_1 e_1'}{2} - \frac{n_1 e_1}{2} + \frac{n_2 e_2'}{2}\right)\right] \qquad (7.5)$$

Squaring both sides of (7.5) and then taking expectations, we get MSE of the estimator $t_{pd}$ up to the first order of approximation as

$$MSE\left(t_{pd}\right) = \overline{Y}^2 f_1\left[C_y^2 + h_1^2 B_1 + h_2^2 B_2 - 2h_1 B_3 - h_2 B_4 + h_1 h_2 B_5\right] \qquad (7.6)$$

Where,

$$\left.\begin{aligned} h_1 &= \frac{4B_2 B_3 - B_4 B_5}{4B_1 B_2 - B_5^2} \\ h_2 &= \frac{2B_1 B_4 - 2B_3 B_5}{4B_1 B_2 - B_5^2} \end{aligned}\right\} \qquad (7.7)$$

and

$$\left.\begin{aligned} B_1 &= f_2 m_2^2 C_{p_2}^2 + f_3 m_1^2 C_{p_1}^2 \\ B_2 &= \frac{1}{4}\left[f_2 n_2^2 C_{p_2}^2 + f_3 n_1^2 C_{p_1}^2\right] \\ B_3 &= f_3 m_1 k_{pb_1} C_{p_1}^2 + f_2 m_2 k_{pb_2} C_{p_2}^2 \\ B_4 &= f_3 n_1 k_{pb_1} C_{p_1}^2 - f_2 n_2 k_{pb_2} C_{p_2}^2 \\ B_5 &= f_3 n_1 m_1 C_{p_1}^2 - f_2 n_2 m_2 C_{p_2}^2 \end{aligned}\right\} \qquad (7.8)$$

## 8. EMPIRICAL STUDY

Data: (Source: Singh and Chaudhary (1986), P. 177).

The population consists of 34 wheat farms in 34 villages in certain region of India. The variables are defined as:

y = area under wheat crop (in acres) during 1974.

$p_1$ = proportion of farms under wheat crop which have more than 500 acres land during 1971. and

$p_2$ = proportion of farms under wheat crop which have more than 100 acres land during 1973.

For this data, we have

N=34, $\overline{Y}$ =199.4, $P_1$ =0.6765, $P_2$ =0.7353, $S_y^2$ =22564.6, $S_{\phi_1}^2$ =0.225490, $S_{\phi_2}^2$ =0.200535,

$\rho_{pb_1}$ =0599, $\rho_{pb_2}$ =0.559, $\rho_\phi$ =0.725.

**Table 7.1: PRE of different estimators of $\overline{Y}$ with respect to $\overline{y}$**

| Choice of scalars | | | | | | | | | |
|---|---|---|---|---|---|---|---|---|---|
| $h_0$ | $h_1$ | $h_2$ | $m_1$ | $m_2$ | $n_1$ | $n_2$ | **Estimator** | **MSE** | **PRE'S** |
| 1 | 0 | 0 | | | | | $\overline{y}$ | 1592.79 | 100 |
| 0 | 1 | 0 | 1 | 0 | | | $t_{d1}$ | 1256.94 | 126.71 |
| | | | 0 | 1 | | | $t_{d2}$ | 1538.00 | 103.90 |
| | | | 1 | 1 | | | $t_{d5}$ | 1197.15 | 133.04 |
| 0 | 0 | 1 | | | 1 | 0 | $t_{d3}$ | 1131.00 | 140.82 |
| | | | | | 0 | 1 | $t_{d4}$ | 2425.83 | 65.65 |
| | | | | | 1 | 1 | $t_{d_6}$ | 1278.00 | 124.62 |
| $h_0$ | $h_1$ | $h_2$ | 1 | 1 | 1 | 1 | $t_{pd}$ | 1032.36 | 154.28 |

## 9. Conclusion

In this paper, we have suggested a class of estimators in single and double sampling by using point bi serial correlation and phi correlation coefficient. From Table 4.1 and Table 7.1, we observe that the proposed estimator $t_p$ and its double sampling version $t_{pd}$, performs better than other estimators considered in this paper.